\documentclass[reqno]{amsart}
\usepackage{pictex}
\usepackage{mathrsfs}
\usepackage{hyperref}
\usepackage{color}
\usepackage{graphicx}
\begin{document}
%\fontsize{11}{4}

%\addtolength{\textwidth}{1in}

%\addtolength{\hoffset}{0in}

%\addtolength{\textheight}{0in}

% \addtolength{\voffset}{0in}

%\setlength{\topmargin}{0in}

%\setlength{\headheight}{0in}

%\setlength{\headsep}{0in}

%\setlength{\textheight}{7.7in}

% \setlength{\textwidth}{7in}

%\setlength{\oddsidemargin}{0in}

% \setlength{\evensidemargin}{0in}

% \setlength{\parindent}{0.25in}

%\setlength{\parskip}{0.25in}

% #############################################
%
%           Macros, etc
%
% #############################################

\def\labelenumi{(\theenumi)}

\newtheorem{thm}{Theorem}[section]
\newtheorem{lem}[thm]{Lemma}
\newtheorem{conj}[thm]{Conjecture}
\newtheorem{cor}[thm]{Corollary}
\newtheorem{add}[thm]{Addendum}
\newtheorem{prop}[thm]{Proposition}
\theoremstyle{definition}
\newtheorem{defn}[thm]{Definition}
\theoremstyle{remark}
\newtheorem{rmk}[thm]{Remark}
\newtheorem{example}[thm]{{\bf Example}}
\newtheorem{Question}[thm]{{ \bf Question}}

\newcommand{\SurfG}{\Sigma_g}
\newcommand{\SurfGN}{\Sigma_{g,n}}
\newcommand{\TriangG}{T_g}
\newcommand{\TriangGOne}{T_{g,1}}
\newcommand{\ProjG}{\mathcal{P}_g}
\newcommand{\TeichG}{\mathcal{T}_g}
\newcommand{\CirclePackGTau}{\mathsf{CPS}_{g,\tau}}
\newcommand{\CrossRatio}{{\bf c}}
\newcommand{\CrossRatioGTau}{\mathcal{C}_{g,\tau}}
\newcommand{\CrossRatioOneTau}{\mathcal{C}_{1,\tau}}
\newcommand{\DeformGTau}{\mathcal{C}_{g,\tau}}
\newcommand{\Forget}{\mathit{forg}}
\newcommand{\Uniform}{\mathit{u}}
\newcommand{\Section}{\mathit{sect}}
\newcommand{\SLTwoC}{\mathrm{SL}(2,{\mathbb C})}
\newcommand{\SLTwoR}{\mathrm{SL}(2,{\mathbb R})}
\newcommand{\SUTwo}{\mathrm{SU}(2)}
\newcommand{\PSLTwoC}{\mathrm{PSL}(2,{\mathbb C})}
\newcommand{\GLTwoZ}{\mathrm{GL}(2,{\mathbb Z})}
\newcommand{\PGLTwoZ}{\mathrm{PGL}(2,{\mathbb Z})}
\newcommand{\GLTwoC}{\mathrm{GL}(2,{\mathbb C})}
\newcommand{\PSLTwoR}{\mathrm{PSL}(2,{\mathbb R})}
\newcommand{\PGLTwoR}{\mathrm{PGL}(2,{\mathbb R})}
\newcommand{\GLTwoR}{\mathrm{GL}(2,{\mathbb R})}
\newcommand{\PSLTwoZ}{\mathrm{PSL}(2,{\mathbb Z})}
\newcommand{\SLTwoZ}{\mathrm{SL}(2,{\mathbb Z})}
\newcommand{\nnn}{\noindent}
\newcommand{\MCG}{{\mathcal {MCG}}}
\newcommand{\MMap}{{\bf \Phi}_{\mu}}
\newcommand{\HH}{{\mathbb H}^2}
\newcommand{\TT}{{T_{1,0}}}
\newcommand{\X}{{\mathcal  X}}
\newcommand{\B}{{\mathcal  B}}
\newcommand{\E}{{\mathcal  E}}
\newcommand{\C}{{\mathscr C}}
\newcommand{\T}{{\mathscr T}}
\newcommand{\M}{{\mathscr M}}
\newcommand{\F}{{\mathscr F}}
\newcommand{\CC}{{\mathbb C}}
\newcommand{\RR}{{\mathbb R}}
\newcommand{\Li}{{{\mathrm{Li}_2}}}

\renewcommand{\L}{{\mathcal  L}}
\newcommand{\G}{{\mathcal  G}}
\newcommand{\R}{{\mathcal  R}}
\newcommand{\Q}{{\mathbb Q}}
\newcommand{\ZZ}{{\mathbb Z}}
\newcommand{\PL}{{\mathscr {PL}}}
\newcommand{\GP}{{\mathcal {GP}}}
\newcommand{\GT}{{\mathcal {GT}}}
\newcommand{\GQ}{{\mathcal {GQ}}}
\newcommand{\EE}{{{\mathcal E}(\rho)}}
\newcommand{\HHH}{{\mathbb H}^3}
\newcommand{\tr}{{\rm tr\, }}

\def\square{\hfill${\vcenter{\vbox{\hrule height.4pt \hbox{\vrule width.4pt
height7pt \kern7pt \vrule width.4pt} \hrule height.4pt}}}$}

\newenvironment{pf}{\noindent {\sl Proof.}\quad}{\square \vskip 12pt}

\title{Solving Thurston Equation in a Commutative Ring}
%for (closed) hyperbolic surfaces]%
%{A new identity for closed hyperbolic surfaces}
\author[F. Luo ]{Feng Luo
 }
\address{Department of Mathematics\\
Rutgers  University\\
Piscataway, NJ 08854, USA}
\email{fluo@math.rutgers.edu}%

%\subjclass[2000]{57M27; 57M50}


\thanks{The work  is partially supported by NSF grants.}

% #############################################
%
%                  Abstract
%
% #############################################

\begin{abstract}
 We show that solutions of Thurston equation on triangulated
3-manifolds in a commutative ring carry topological information.
We also introduce a homogeneous Thurston equation and a
commutative ring associated to triangulated 3-manifolds.
\end{abstract}

\maketitle

\section{introduction}

\subsection{Statement of  results}
Given a triangulated oriented 3-manifold (or pseudo 3-manifold)
$(M, \mathcal T$), Thurston equation associated to $\mathcal T$ is
a system of integer coefficient polynomials. W. Thurston \cite{Th}
introduced his equation in the field $\mathbf C$ of complex
numbers in order to find hyperbolic structures. Since then, there
have been much research on solving Thurston equation in $\mathbf
C$, \cite{NZ}, \cite{Cho}, \cite{PW}, \cite{Ti}, \cite{Yo},
\cite{BB} and others.  Since the equations are integer coefficient
polynomials, one could attempt to solve Thurston equation in a
ring with identity. The purpose of this paper is to show that
interesting topological results about the 3-manifolds can be
obtained by solving Thurston equation in a commutative ring with
identity. For instance, by solving Thurston equation in the field
$\mathbf Z/\mathbf 3Z$ of three elements, one obtains the result
which was known to H. Rubinstein and S. Tillmann that a closed
1-vertex triangulated 3-manifold is not simply connected if each
edge has even degree.

\bigskip
\begin{thm}\label{thm:main} Suppose $(M, \mathcal T)$ is an
oriented connected closed 3-manifold with a triangulation
$\mathcal T$ and $R$ is a commutative ring with identity. If
Thurston equation on $(M, \mathcal T)$ is solvable in $R$ and
$\mathcal T$ contains an edge which is a loop, then there exists a
homomorphism from $\pi_1(M)$ to $PSL(2, R)$ sending the loop to a
non-identity element. In particular, $M$ is not simply connected.
%\begin{equation}\label{eqn:maintheorem}
   %\end{equation}
\end{thm}

We remark that the existence of an edge which is a loop  cannot be
dropped in the theorem. Indeed, it was observed in  \cite{KK},
 \cite{Ti}, and \cite{Yang} that for simplicial triangulations
 $\mathcal T$ and any commutative ring $R$, there are always
 solutions to Thurston equation.  Theorem 1.1 for $R =\mathbf C$
 was first proved by Segerman-Tillmann \cite{ST}. A careful
 examination of the proof of \cite{ST} shows that their method
 also works for any field $R$. However, for a commutative ring
 with zero divisors, the geometric argument breaks down.
 We prove theorem 1.1 by introducing a homogeneous
 Thurston equation and studying its solutions.
 Theorem 1.1 prompts us to introduce the universal construction of a
 Thurston ring of a triangulated 3-manifold.
 Theorem 1.1 can be phrased in terms of the universal
 construction (see theorem 5.2).

Thurston equation can be defined for any ring (not necessary
commutative) with identity (see \S2). We do not know if theorem
1.1 holds in this case. The most interesting non-commutative rings
for 3-manifolds are probably the algebras of $2 \times 2$ matrices
with real or complex coefficients. Solving Thurston equation in
the algebra  $M_{2 \times 2}(\mathbf R)$ has the advantage of
linking hyperbolic geometry to $Ad(S^3)$ geometry. See \cite{Da}
for related topics.

Motivated by theorem 1.1, we propose the following two
conjectures.

\medskip

\noindent {\bf Conjecture 1.} If $M$ is a compact 3-manifold and
$\gamma \in \pi_1(M) -\{1\}$, there exists a finite commutative
ring $R$ with identity and a homomorphism from $\pi_1(M)$ to
$PSL(2, R)$ sending $\gamma$ to a non-identity element.

\medskip

 \noindent {\bf Conjecture 2.} If $M \neq S^3$ is a closed
oriented 3-manifold, then there exists a 1-vertex triangulation
$\mathcal T$ of $M$ and a commutative ring $R$ with identity so
that Thurston equation associated to $\mathcal T$ is solvable in
$R$.

\medskip

Conjecture 2 is supported by the main result in \cite{LTY}. It
states that if $M$ is closed hyperbolic and $\mathcal T$ is a
1-vertex triangulation so that all edges are homotopically
essential, then Thurston's equation on $\mathcal T$ is solvable in
$\mathbf C$.

\subsection{Organization of the paper}
In \S2, we recall briefly triangulations of 3-manifolds and pseudo
3-manifolds and Thurston equation. A homogeneous Thurston equation
is introduced. In \S3, we recall cross ratio and projective plane
in a commutative ring. Theorem 1.1 is proved in \S4. In \S5, we
introduce a universal construction of Thurston ring of a
triangulated 3-manifold and investigate the relationship between
Thurston ring and Pachner moves. Some examples of solutions of
Thurston equation in finite rings are worked out in \S6.

\subsection{Acknowledgement}
We thank  Sergei Matveev and Carlo Petronio for answering my
questions on triangulations of 3-manifolds.

%\begin{figure}[ht!]
%\centering
%\includegraphics[scale=0.55]{2.eps}
%\caption{3-holed spheres and 1-holed tori} \label{figure 1}
%\end{figure}

%\subsection{Basic idea of the proof}

%\end{document}

\section{Preliminaries on Triangulations and Thurston Equation}

All manifolds and tetrahedra are assumed to be oriented in this
paper. We assume all rings have the identity element.

\subsection{Triangulations}
A compact oriented triangulated pseudo 3-manifold $(M, \mathcal
T)$ consists of a disjoint union $X =\sqcup_{i} \sigma_i$ of
oriented Euclidean tetrahedra $\sigma_i$ and a collection of
orientation reversing affine homeomorphisms $\Phi$ between some
pairs of codimension-1 faces in $X$. The pseudo 3-manifold $M$ is
the quotient space $X/\Phi$ and the simplices in $\mathcal T$ are
the quotients of simplices in $X$. The \it boundary \rm $\partial
M$ of $M$ is the quotient of the union of unidentified
codimension-1 faces in $X$. If $\partial M =\emptyset$, we call
$M$ \it closed\rm.  The sets of all quadrilateral types (to be
called \it quads \rm for simplicity) and normal triangle types in
$\mathcal T$ will be denoted by $\Box = \Box(\mathcal T)$ and
$\triangle =\triangle(\mathcal T)$ respectively. See \cite{JR},
\cite{Ha} or \cite{Lu1} for more details. The most important
combinatorics ingredient in defining Thurston equation is a
$\mathbf Z/ 3\mathbf Z$ action on $\Box$ which we recall now. If
edges of an oriented tetrahedron $\sigma$ are labelled by $a,b,c$
so that opposite edges are labelled by the same letters (see
figure 1(a)), then the cyclic order $a \to b \to c \to a$ viewed
at vertices is independent of the choice of the outward normal
vectors at vertices and depends only on the orientation of
$\sigma$. Since a quad in $\sigma$ corresponds to a pair of
opposite edges, this shows that there is a $\mathbf Z /3 \mathbf
Z$ action on the set of all quads in $\sigma$ by cyclic
permutations. If $q, q' \in \Box$, we use $q \to q'$ to indicate
that $q,q'$ are in the same tetrahedron so that $q$ is ahead of
$q'$ in the cyclic order. The set of all i-simplices in $\mathcal
T$ will be denoted by $\mathcal T^{(i)}$. Given an edge $e \in
\mathcal T^{(1)}$, a tetrahedron  $\sigma \in \mathcal T^{(3)}$
and a quad $q \in \Box(\mathcal T)$, we use $e < \sigma$ to
indicate that $e$ is adjacent to $\sigma$; $q \subset \sigma$ to
indicate that $q$ is inside $\sigma$; and $q \sim e$ to indicate
that $q$ and $e$ are in the same tetrahedron and $q$ faces $e$ (in
the unidentified space $X$). See \cite{Lu1} for more details. An
edge $e \in \mathcal T^{(1)}$ is called \it interior \rm if it is
not in the boundary $\partial M$. In particular, if $(M, \mathcal
T)$ is a closed triangulated pseudo 3-manifold, all edges are
interior.

\subsection{Thurston equation}
\begin{defn}\label{thurstonequ}
Give a compact triangulated oriented pseudo 3-manifold $(M,
\mathcal T)$ and a ring $R$ (not necessary commutative) with
identity, a function $x:\Box(\mathcal T) \to R$ is called a
solution to Thurston equation associated to $\mathcal T$ if

(1) whenever $q \to q'$ in $\Box$, $x(q')(1-x(q))=(1-x(q))x(q')
=1$,

(2) for each interior edge $e \in \mathcal T^{(1)}$ so that $q_1,
..., q_n$ are quads facing $e$ labelled cyclically around $e$,
$$ x(q_1).....x(q_n) =1,  \quad  x(q_n) ... x(q_1) =1.$$
\end{defn}

%\begin{prop}\label{prop:loopisembedded}
%\end{prop}

\begin{figure}[ht!]
\centering
\includegraphics[scale=0.6]{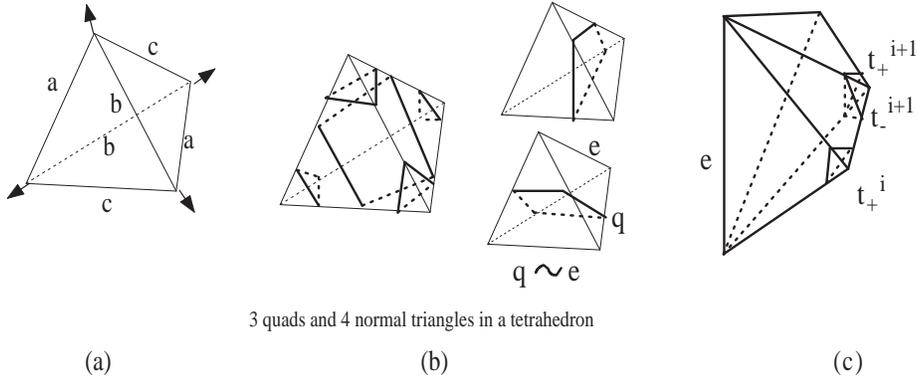}
\caption{cyclic order on three quads in a tetrahedron}
\label{figure 1}
\end{figure}

Note that condition (1) implies both $x(q)$ and $x(q)-1$ are
invertible elements in $R$ with inverses $1-x(q'')$ and $x(q')$
where $q'' \to q$. If the ring $R$ is commutative which will be
assumed from now on, we only need one equation in each of (1) and
(2).

\begin{example}  If $R =\mathbf Z/3 \mathbf Z =\{0,1,2\}$ is the
field of 3 elements. Then a solution $x$ to Thurston equation must
satisfy $x(q) =2$ for all $q$.  In this case, the first condition
(1) holds. The second equation at an edge $e$ becomes $2^k=1$
where $k$ is the degree of  $e$. Since $2^k=1$ if and only if $k$
is even, we conclude that Thurston equation has a solution in
$\mathbf Z/3 \mathbf Z$ if and only if each interior edge has even
degree.
\end{example}

A related homogeneous version of the equation is

\begin{defn} (Homogeneous Thurston Equation (HTE)) Suppose $(M, \mathcal T)$ is a
 compact oriented triangulated
pseudo 3-manifold and  $R$ is a commutative ring with identity. A
function $z: \Box \to R$ is called a solution to the homogeneous
Thurston equation if

(1) for each tetrahedron $\sigma \in \mathcal T$, $\sum_{ q
\subset \sigma} z(q) =0$,

(2) for each interior edge $e$ so that the set of all quads facing
it is $\{q_1,..., q_n\}$,
$$ \prod_{i=1}^n z(q_i) = \prod_{i=1}^n (-z(q_i'))$$ or simply
$$  \prod_{ q \sim e} z(q) = \prod_{ q \sim e} (-z(q')),$$ where
$q \to q'$.
\end{defn}

Note that if $z$ solves HTE and $k: \mathcal T^{(3)} \to R$ is any
function, then $w(q) = k(\sigma) z(q): \Box \to R$, $q \subset
\sigma$, is another solution to HTE.  Let $R^*$ be the group of
all invertible elements in $R$. If $z$ solves HTE and $z(q) \in
R^*$ for all $q \in \Box$, then $x(q) = -z(q) z(q')^{-1}$ for $q
\to q'$ solves Thurston equation. Indeed, condition (2) in
definition 2.1 follows immediately from condition (2) in
definition 2.3 by division. To check condition (1), suppose $q \to
q' \to q'' \to q$. Then $z(q'') = -z(q) - z(q')$. Furthermore,
$x(q) = -z(q)/z(q')$ and $x(q') = -z(q')/z(q'') = z(q')/(z(q) +
z(q')) = 1/(1-x(q)).$  Conversely, we have,

\begin{lem} If $R$ is a commutative ring with identity and $x:\Box
\to R$ solves Thurston equation, then there exists a solution $z:
\Box \to R^*$ to HTE so that for all $q \in \Box$, $x(q) =
-z(q)/z(q')$. \end{lem}

\begin{proof}  By definition, each $x(q)$ is invertible. For each tetrahedron $\sigma$ containing three
quads $q_1 \to q_2 \to q_3$, we have $x(q_2) = 1/(1-x(q_1)),
x(q_3) =(x(q_1)-1)/x(q_1)$. Define a map $z: \Box \to R^*$ by
$z(q_1)= x(q_1), z(q_2) = -1, z(q_3) = 1-x(q_1)$. Then by
definition $x(q) = -z(q)/z(q')$ for all $q$ and $\sum_{ q \subset
\sigma} z(q) =0$ for each $\sigma \in \mathcal T^{(3)}$. Due to
$x(q)=-z(q)/z(q')$ and $\prod_{ q \sim e} x(q) =1$, we see that
$\prod_{ q \sim e} z(q) = \prod_{q \sim e} (-z(q'))$.  \end{proof}

We remark that the solution $z$ in the lemma depends on the
specific choice of the quad $q_1$ in each tetrahedron.

There is a similar version of homogeneous Thurston equation where
we replace the condition $  \prod_{ q \sim e} z(q) = \prod_{ q
\sim e} (-z(q'))$ at interior edge $e$ by $  \prod_{ q \sim e}
z(q'') = \prod_{ q \sim e} (-z(q'))$ where $q'' \to q \to q'$. In
this setting, the transition from HTE to Thurston equation is
given by $x(q) = -z(q')/z(q'')$.
%\end{document}

\section{Cross ratio and projective line in a commutative ring}

Let $R$ be a commutative ring with identity and $R^*$ be the group
of invertible elements in $R$. Let $GL(2,R)$ and
$PGL(2,R)=GL(2,R)/\sim$ where $M \sim \lambda M$, $\lambda \in
R^*$ be the general linear group and its projective group. The
group $GL(2,R)$ acts linearly from the left on $R^2 =\{ \left(
\begin{array}{c}
a\\
b\\
\end{array}
\right) |a,b, \in R\}$. Define the skew symmetric bilinear form
$<,>$ on $R^2$ by $<\left(
\begin{array}{c}
a\\
b\\
\end{array} \right),
 \left(
\begin{array}{c}
c\\
d\\
\end{array} \right)> = ad-bc$ which is the determinant
of $\left[\begin{array}{cc} a & c \\ b & d \end{array}\right].$
Note that for $X= \left[\begin{array}{cc} a & c \\ b & d
\end{array}\right]$, its adjacent matrix
$adj(X)=\left[\begin{array}{cc} d & - c
 \\- b & a \end{array}\right] $ satisfies $X adj(X) = det(X) I$.
 We also use the transpose to write $\left(\begin{array}{c}
a\\
b\\
\end{array} \right)$ as $(a,b)^t$.

By the basic properties of the determinant, we have
\begin{lem} Suppose $A,B, A_1, ..., A_n \in R^2$ and $X \in
GL(2,R)$. Then

(1) $< A, B> =-<B,A>$ and $<XA, XB> =det(X)<A,B>$;

(2) If $A =\left(\begin{array}{c} a\\ b\\ \end{array} \right)$ and
$B =\left(\begin{array}{c} c\\ d\\ \end{array} \right)$ so that
$<A, B> \in R^*$, then $XA =\left(\begin{array}{c} 1\\ 0\\
\end{array} \right)$ and $XB =\left(\begin{array}{c} 0\\ 1\\ \end{array}
\right)$ where $X =\frac{1}{<A,B>}\left[\begin{array}{cc} d & -c\\
-b & a
\end{array}\right]$;

(3) $<A_1, A_2> A_3 + <A_2, A_3>A_1 + < A_3, A_1>A_2=0$;

(4) Let $R_{ijkl} =< A_i, A_j><A_k, A_l>$. Then
$R_{ijkl}=R_{jilk}=R_{klij} = -R_{jikl} = -R_{ijlk}$ and
$R_{ijkl}+R_{iklj}+R_{iljk}=0$.
\end{lem}

Indeed, (1) and (2) follow from the properties of determinant. To
see (3), let $A_i=\left(\begin{array}{c} a_i\\ b_i\\ \end{array}
\right)$, then the left-hand-side of (3) is of the form
$\left(\begin{array}{c} x\\ y\\
\end{array}\right)$ where $x =det \left[\begin{array}{ccc}
a_1&a_2&a_3 \\ b_1&b_2&b_3 \\ a_1 &a_2 &a_3 \end{array}\right]$
and $y=det \left[\begin{array}{ccc} a_1&a_2&a_3 \\ b_1&b_2&b_3 \\
b_1&b_2&b_3 \end{array}\right]$ by the row expansion formula for
determinants. Now these two $3 \times 3$ determinants are zero.
Thus (3) follows. The first set of identities for $R_{ijkl}$
follow from the definition. The second follows from (3) by
applying the bilinear form $<,>$ to it with $A_i$. We remark that
$R_{ijkl}$ enjoys the same symmetries that a Riemannian curvature
tensor does.

\begin{defn}(Cross ratio) Suppose $A_1, ..., A_4 \in R^2$. Then
their cross ratio, denoted by $(A_1, A_2; A_3, A_4)$ is defined to
be the vector $\left(\begin{array}{c} R_{1423}\\ R_{1324}\\
\end{array} \right)   =  \left(\begin{array}{c} R_{1423}\\ -R_{1342}\\
\end{array} \right) \in R^2   $ where $R_{ijkl} =<A_i, A_j><A_k,
A_l>$. \end{defn}
For instance, \begin{equation} <\left(\begin{array}{c} 1\\ 0\\
\end{array} \right), \left(\begin{array}{c} 0\\ 1\\ \end{array}
\right); \left(\begin{array}{c} a\\ b\\ \end{array} \right),
\left(\begin{array}{c} x\\ y\\ \end{array} \right)>
=\left(\begin{array}{c} -ay\\ -bx\\ \end{array} \right)
\end{equation}

By lemma 3.1, we obtain,

\begin{cor} Suppose $A_1, ..., A_n \in R^2$. Then

(1) $(A_1, A_2;A_3, A_4) =(A_3, A_4; A_1, A_2) =(A_2, A_1; A_4,
A_3)$,

(2) if $(A_1, A_2; A_3, A_4) =\left(\begin{array}{c} a\\ b\\
\end{array} \right)$, then $(A_2, A_1; A_3, A_4)$
$= \left(\begin{array}{c} b\\ a\\ \end{array} \right)$,

(3) $(A_1, A_2; A_3, A_4) + (A_1, A_3; A_4, A_2)+ (A_1, A_4; A_2,
A_3) =0$,

(4) if $X \in GL(2,R)$, then $(XA_1, XA_2; XA_3, XA_4) = det(X)^2
(A_1, A_2; A_3, A_4)$,

(5) if $B,C \in R^2$ and $A_{n+1}=A_1$, then
$$\prod_{i=1}^n (<B,A_{i+1}><C, A_i> ) = \prod_{i=1}^n(<B,A_i><C,
A_{i+1}>).$$ \end{cor}

Corollary 3.3 shows that the cross ratio $(A_1, A_2; A_3, A_4)
=(A_i, A_j; A_k, A_l)$ whenever $\{i,j\}=\{1,2\}, \{k,j\}=\{3,4\}$
and $(i,j,k,l)$ is a positive permutation of $(1,2,3,4)$, i.e.,
the cross ratio depends only on the partition $\{i,j\} \sqcup
\{k,l\}$ of $\{1,2,3,4\}$ and the orientation of $(1,2,3,4)$. This
shows that if $\sigma$ is an oriented tetrahedron so that its i-th
vertex is assigned a vector $A_i \in R^2$, then one can define the
cross ratio of a quad $q \subset \sigma$ to be $(A_i, A_j; A_k,
A_l)$ where $q$ corresponds to the partition $\{i,j\} \sqcup
\{k,l\}$ of the vertex set $\{1,2,3,4\}$ and $(i,j,k,l)$
determines the orientation of $\sigma$.

\begin{example} (Solutions of HTE by cross ratio) Given any compact triangulated
pseudo 3-manifold $(M, \mathcal T)$ and $f: \triangle \to R^2$ so
that $f(t) = f(t')$ when two normal triangles $t,t'$ share a
common normal arc, we define a map $F: \Box(\mathcal T) \to R^2$
by $F(q) =(f(t_1),
f(t_2); f(t_3), f(t_4)) = \left(\begin{array}{c} z(q)\\ y(q)\\
\end{array} \right)$ where $t_1, .., t_4$ are the four normal
triangles in a tetrahedron $\sigma$ containing $q$ so that $q$
separates $\{t_1, t_2\}$ from $\{t_3, t_4\}$ and $t_1 \to t_2 \to
t_3 \to t_4$ defines the orientation of $\sigma$. Then corollary
3.3 shows that $z : \Box \to R$ is a solution to HTE. Note that
$y(q)=-z(q')$ with $q \to q'$.
\end{example}

This example serves as a guide for us to solve Thurston equation
and HTE. Indeed, the goal is to solve Thurston equation by writing
each solution $x : \Box \to R$ in terms of cross ratio in the
universal cover $(\tilde{M}, \tilde{ \mathcal T})$.

Let $R^2/\sim$ be the quotient space where $u \sim \lambda u$, $u
\in R^2$ and $\lambda \in R^*$. If $ x =(a,b)^t \in R^2$, then
$[x]=[a,b]^t$ denotes the image of $x$ in $R^2/\sim$.

\begin{defn} The \it projective line \rm $PR^1 = \{A \in
R^2 | \text{there exists $B \in R^2$ so that} <A, B> \in
R^*\}/\sim$ where $A \sim \lambda A$ for $\lambda \in R^*$.  A set
of elements $\{A_1, ..., A_n\} $ (or \{$[A_1], ..., [A_n]$\}) in
$R^2$ (or in $PR^1$) is called \it admissible \rm if $<A_i, A_j>
\in R^*$ for all $i \neq j$. The \it cross ratio \rm of four
points $\alpha_i, i=1,2,3,4$ in $PR^1$, denoted by $[\alpha_1,
\alpha_2; \alpha_3, \alpha_4]$, is the element $[(A_1, A_2;A_3,
A_4)] \in R^2/\sim$ so that $A_i \in \alpha_i$. We also use $[A_1,
A_2; A_3, A_4] \in PR^1$ to denote $[(A_1, A_2;A_3, A_4)]$.
\end{defn}

\begin{prop} (1) Given an admissible set of three elements $A_1,
A_2, A_3 \in R^2$ and $v =\left(\begin{array}{c} c_1\\ c_2\\
\end{array} \right)$, there exists a unique $A_4 \in R^2$ so that
$(A_1, A_2; A_3, A_4) =v$. Furthermore, $A_1, ..., A_4$ form an
admissible set if and only if $c_1, c_2, c_1-c_2 \in R^*$.

(2) Suppose $A_1, ..., A_4, B_1, ..., B_4 \in R^2$ so that both
$\{A_1, ..., A_4\}$ and $\{B_1, ..., B_4\}$ are admissible and
$[A_1, A_2; A_3, A_4]=[B_1, B_2; B_3, B_2]$. Then there exists a
unique $X \in PGL(2,R)$ so that $[XA_i ]=[B_i]$ for all $i$.
Furthermore, if $Y \in GL(2,R)$ so that $[YA_i]=[A_i], i=1,2,3$,
then $Y = \lambda I$ for $\lambda \in R^*$.
\end{prop}

\begin{proof} To see the existence part of (1), let $A_i =(a_i,
b_i)^t$
% \left(\begin{array}{c} a_i\\ b_i\\ \end{array} \right)$
and consider $X =\frac{1}{<A_1, A_2>}$ $ \left(\begin{array}{cc}
b_2&-a_2\\-b_1&a_1 \end{array}\right)$ $ \in GL(2,R)$. Then
$XA_1=(1,0)^t$  and $XA_2=(0,1)^t$. By corollary 3.3(4), after
replacing $A_i$ by $XA_i$, we may assume that $A_1=(1,0)^t$,
$A_2=(0,1)^t$.  Then by identity (1), $(A_1, A_2; A_3, A_4)
=(-a_3b_4, -a_4b_3)^t$. By the assumption that $<A_i, A_3> \in
R^*$ for $i=1,2$, we see that $a_3, b_3 \in R^*$. It follows that
$a_4=-c_1/b_3$ and $b_4 =-c_2/a_3$. This shows that $A_4$ exists
and is unique.

Now given that $\{A_1, A_2, A_3\}$ is admissible, the set $\{A_1,
..., A_4\}$ is admissible if and only if $<A_4, A_i> \in R^*$,
i.e., $c_1, c_2, c_1-c_2 \in R^*$.

To see part (2), by the proof of part (1) and the assumption
$<A_1, A_2>, <B_1, B_2> \in R^*$, after replacing $A_i$ by $YA_i$
and $B_i$ by $ZB_i$ for some $Y,Z \in GL(2,R)$, we may assume that
$A_1=B_1=(1,0)^t$ and $A_2=B_2 =(0,1)^t$. Let $A_3=(a,b)^t,
A_4=(a',b')^t$, $B_3 =(c,d)^t$ and $B_4=(c',d')^t$. Then the
admissible condition implies that $a,b,c,d,a',b',c',d' \in R^*$.
Furthermore, $[A_1, A_2;A_3, A_4] =[B_1, B_2;B_3, B_4]$ implies
that there exists $\lambda \in R^*$ so that $a' =\frac{\lambda d
d'}{b}$ and $b' =\frac{\lambda c c'}{a}$. This
shows that the matrix $X=\left[\begin{array}{cc} \frac{c}{a}&0\\
0&\frac{d}{b}
\end{array}\right] \in GL(2,R)$ satisfies $XA_1 = \frac{c}{a}B_1$,
$XA_2 = \frac{d}{b} B_2$, $XA_3=B_3$ and $XB_4 = \lambda
\frac{cd}{ab}B_4$.

To see the uniqueness, say $YA_i = \lambda_i A_i$ for $\lambda_i
\in R^*$. We claim that $\lambda_1 =\lambda_2 =\lambda_3$ and
$Y=\lambda_1 I$. Indeed, by definition, $det(Y) <A_i, A_j> =<
YA_i, YA_j> =\lambda_i \lambda_j < A_i, A_j>$. Since $<A_i, A_j>
\in R^*$, we obtain $\lambda_i \lambda_j =det (Y)$. This implies
that $\lambda_i =\lambda_1$ for $i=2,3$. We conclude that
$Y[A_1,A_2] =\lambda_1[A_1, A_2]$. Since the matrix $[A_1, A_2]
\in GL(2,R)$, it follows that $Y =\lambda_1 I$.

\end{proof}

\section{A proof of theorem 1.1}

We will prove a slightly general theorem which holds for compact
oriented pseudo 3-manifolds $(M, \mathcal T)$. Let $M^*$ be $M$
with a small regular neighborhood of each vertex removed and let
$\mathcal T^*$ be the ideal triangulation $\{s \cap M^* | s \in
\mathcal T\}$ of the compact 3-manifold $M^*$.

\begin{thm} Suppose $(M, \mathcal T)$ is a compact triangulated
pseudo 3-manifold and $R$ is a commutative ring with identity so
that Thurston equation on $\mathcal T$ is solvable in $R$. Then
each edge $e \in \mathcal T^*$ lifts to an arc in the universal
cover $\tilde{M}^*$ of $M^*$ joining different boundary components
of $\tilde{M}^*$. Furthermore, if $M$ is a closed  connected
3-manifold so that there exists an edge $e$ having the same end points, then there exists a
representation of $\pi_1(M)$ into $PSL(2,R)$ sending the loop $[e]$
to a non-identity element. \end{thm}

The main idea of the proof is based on the methods developed in
\cite{NZ}, \cite{ST}, \cite{Ti},  and \cite{Yo} which construct
pseudo developing map and the holonomy associated to a solution to
Thurston equation.

\subsection{Pseudo developing map}
Let $\pi: \tilde{M}^* \to M^*$ be the universal cover and
$\tilde{\mathcal T^*}$ be the pull back of the ideal triangulation
$\mathcal T^*$ of $M^*$ to $\tilde{M}^*$. We use
$\tilde{\triangle}$ and $\tilde{\Box}$ to denote the sets of all
normal triangle types and quads in $\tilde{\mathcal T^*}$
respectively. The sets of all normal triangle types and quads in
$\mathcal T^*$ are the same as those of $\mathcal T$ and will
still be denoted by $\triangle$ and $\Box$. The covering map $\pi$
induces a surjection $\pi_*$ from $\tilde{\triangle}$ and
$\tilde{\Box}$ to $\triangle$ and $\Box$ respectively so that
$\pi_*(d_1)=\pi_*(d_2)$ if and only if $d_1$ and $d_2$ differ by a
deck transformation element.

Suppose $x : \Box \to R$ solves Thurston equation on $\mathcal T$
and $z: \Box \to R^*$ is an associated solution to HTE constructed
by lemma 2.1. Let $w: \Box \to PR^1$ be the map $w(q) =[z(q),
-z(q')]^t$ where $q \to q'$. Let $\tilde{x} = x \pi_*$, $\tilde{z} =
z \pi_*$ and $\tilde{w}= w \pi_*$ be the associated maps defined on
$\tilde{\Box}$. By the construction, $\tilde{x}$ and $\tilde{z}$
are solutions to Thurston equation and HTE on $\tilde{\mathcal
T^*}$.

\begin{defn} (See \cite{Ti}, \cite{ST}, \cite{Yo}, \cite{NZ})
Given a solution  $x$ to Thurston equation on $(M, \mathcal T)$, a
map $\phi: \tilde{\triangle} \to PR^1$ is called a \it pseudo
developing map \rm associated to $x$ if

(1) whenever $t_1, t_2$ are two normal triangles in
$\tilde{\triangle}$ sharing a normal arc, denoted by $t_1 \sim
t_2$ in the sequel, then $\phi(t_1) =\phi(t_2)$,

(2) if $t_1, t_2, t_3, t_4$ are four normal triangles in a
tetrahedron $\sigma$ then \{$\phi(t_1), .$$.., \phi(t_4)$ \} is
admissible and  \begin{equation} [\phi(t_1), \phi(t_2); \phi(t_3),
\phi(t_4)] = \tilde{w}(q) \end{equation} where $t_1 \to t_2 \to
t_3 \to t_4$ determines the orientation of the tetrahedron
$\sigma$ and $q \subset \sigma$
 is the quad separating $\{t_1, t_2\}$ from $\{t_3, t_4\}$.
 \end{defn}

The main result in this section is

\begin{thm} Given any solution $x$ to Thurston equation on a
compact pseudo 3-manifold $(M, \mathcal T)$, there exists a pseudo
developing map associated to $x$. \end{thm}

\begin{proof}
The proof is based on the following result which is an immediate
consequence of proposition 3.6.

\begin{lem} Suppose $\{t_1, t_2, t_3, t_4\}$ are four normal
triangles in $\sigma \in \tilde{\mathcal T}^{(3)}$ so that $t_1
\to t_2 \to t_3 \to t_4$ determines the orientation and $q$ is a
quad in $\sigma$ separating $\{t_1, t_2\}$ from $\{t_3, t_4\}$. If
$\phi(t_i) \in PR^1$, $i=1,2,3$, are defined and
$\{\phi(t_1),\phi(t_2), \phi(t_3)\}$ is admissible, then there
exists a unique $\phi(t_4) \in PR^1$ so that $[\phi(t_1),
\phi(t_2); \phi(t_3), \phi(t_4)]$ $ = \tilde{w}(q)$. Furthermore,
$\{\phi(t_1), ..., \phi(t_4)\}$ is admissible. \end{lem}

Indeed, the existence and uniqueness follows from proposition 3.6.
The admissibility of $\{\phi(t_1), ..., \phi(t_4)\}$ follows from
that fact that $z(q), -z(q'), z(q)-(-z(q'))=-z(q'')$ are in $R^*$
where $\tilde{w}(q) =[z(q), -z(q')]^t$ and $q \to q' \to q''$.

We now use the lemma to construct the pseudo developing map $\phi$
by the ``combinatorial continuation" method. To begin, by working
on connected component of $M$, we may assume that $M$ is connected.
 Let $G$ be the connected graph dual to
the ideal triangulation $\tilde{\mathcal T^*}$ of $\tilde{M}^*$,
i.e., vertices of $G$ are tetrahedra in $\tilde{\mathcal T^*}$ and
edges in $G$ are pairs of tetrahedra sharing a codimension-1 face.
An edge path $\alpha =[\sigma_1, ..., \sigma_n; \tau_1, ...,
\tau_{n-1}]$ in $G$ consists of tetrahedra $\sigma_i$ and
codimension-1 faces $\tau_i$ so that $\tau_i \subset \sigma_i \cap
\sigma_{i+1}$. If $\sigma_{n}=\sigma_1$, we say $\alpha$ is an
edge loop.

\begin{lem} Suppose $\alpha$ is an edge path from $\sigma_1$ to
$\sigma_n$ and $t_1, t_2, t_3$ are normal triangles in $\sigma_1$
adjacent to the codimension-1 face $\tau_1$ so that $\phi(t_1),
\phi(t_2), \phi(t_3)$ are defined and admissible. Then there
exists an extension $\phi_{\alpha}$ of $\phi$ (depending on
$\alpha$) to all normal triangles $t$ in $\sigma_i$'s so that
conditions (1) and (2) in definition 4.2 hold.
\end{lem}

\begin{proof}  Suppose $t_4$ is the last normal triangle in
$\sigma_1$. By lemma 4.4, we define $\phi_{\alpha}(t_4)$ so that
identity (2) holds for the quad $q$ separating $t_1, t_2$ from
$t_3, t_4$ (subject to orientation). Note that this implies that
identity (2) holds for all other quads $q^*$ in $\sigma_1$ due to
the basic property of cross ratio (corollary 3.3) and that $x$
solves Thurston equation.  Now we extend $\phi_{\alpha}$ to normal
triangles in $\sigma_2$ as follows. Suppose $t_1', t_2', t_3'$ are
the normal triangles in $\sigma_2$ so that $t_i' \sim t_i$, i.e., they share
a normal arc. Define $\phi_{\alpha}(t_i') =\phi_{\alpha}(t_i)$ and
then use lemma 4.4 to extend $\phi_{\alpha}$ to the last normal triangle in
$\sigma_2$. Inductively, we define $\phi_{\alpha}$ for all normal
triangles $t$ in $\sigma_i$. By the construction, both conditions
(1) and (2) in definition 4.2 hold for $\phi_{\alpha}$.
\end{proof}

We call $\phi_{\alpha}$ the ``combinatorial continuation" of $\phi$
along the edge path $\alpha$ and denote it by
$\phi^{\sigma_1}_{\alpha}$ to indicate the initial value. From the
construction, if $\beta$ is an edge path starting from $\sigma_n$
to $\sigma_m$ and $\beta \alpha$ is the multiplication of the edge
paths $\alpha$ and $\beta$, then
\begin{equation}
\phi^{\sigma_n}_{\beta} (t) = \phi^{\sigma_1}_{\beta \alpha}(t)
\end{equation}
for all normal triangles $t$ in $\sigma_m$.

Our goal is to show that the extension $\phi_{\alpha}^{\sigma_1}$
is independent of the choice of edge path $\alpha$, i.e.,
$\phi_{\alpha}^{\sigma_1}(t)=\phi_{\alpha'}^{\sigma_1}(t)$ for two
edge paths $\alpha$ and $\alpha'$ from $\sigma_1$ to $\sigma_n$
and $t \subset \sigma_n$. By (3), this is the same as showing
$\phi_{\alpha' \alpha^{-1}}^{\sigma_n}(t)=t$. Therefore,
  it suffices to show

\begin{lem} If $\alpha$ is an edge loop in $G$ from
$\sigma_1$ to $\sigma_1$, then $\phi^{\sigma_1}_{\alpha}(t) = t$
for all normal triangles $t$ in $\sigma_1$.
\end{lem}
\begin{proof} Form the 2-dimensional CW complex $W$ by attaching
2-cells to the graph $G$ as follows. Recall that an edge $e \in
\tilde{\mathcal T^*}^{(3)}$ is called \it interior \rm if it is
not in the boundary $\partial \tilde{M}^*$. For each interior edge
$e$ in $\tilde{\mathcal T^*}$ adjacent to tetrahedra $\delta_1,
..., \delta_m$, ordered cyclically around $e$, there corresponds
an edge loop $\alpha_e =[\delta_1, ..., \delta_m; \epsilon_1, ...,
\epsilon_m]$ where $\epsilon_i \subset \delta_i \cap \delta_{i+1}$
and $\delta_{m+1} =\delta_1$. We attach a 2-cell to $G$ along
$\alpha_e$ for each interior edge $e$ to obtain $W$. By the
construction, the universal cover space $\tilde{M}^*$ is obtained
from $W$ by attaching a product space $B \times [0, 1)$ along a
surface $B \times 0$. Thus $W$ is homotopic to $\tilde{M}^*$. In
particular, $W$ is simply connected. This shows that the edge loop
$ \alpha$ is a product of edge loops of the form $\alpha_e$, for
interior edges $e$, and loops of the form $\beta \beta^{-1}$ for
some edge paths $\beta$. By the identity (3), the lemma holds for
edge loops of the form $\beta \beta^{-1}$. Therefore, it remains
to prove the lemma for edge loops $\alpha = \alpha_e = [\delta_1,
..., \delta_m; \epsilon_1, ..., \epsilon_m]$.

To this end, suppose that $\phi$ is defined at the normal
triangles $t^1_0, t^1_+, t^1_{\infty}$ in the tetrahedron $\delta_1$
so that the edge $e$ is adjacent to $t^1_0, t^1_{\infty}$. Let the
normal triangles in the tetrahedron $\delta_i$ be $t_0^i, t_+^i,
t_-^{t+1}, t^i_{\infty}$ so that $t^i_0, t^i_{\infty}$ are
adjacent to the edge $e$ and $t^i_0 \to t^i_{\infty} \to t_+^i \to
t^{i+1}_-$ defines the orientation. Then by the construction,
 $t_0^i \sim t_0^{i+1}$, $t_{\infty}^i \sim
t_{\infty}^{i+1}$ and $t^i_{+} \sim t^i_{-}$ where indices are
counted modulo $m$. See figure 1(c). Let $q_i$ be the quad in
$\delta_i$ separating $\{t_0^i, t_{\infty}^i\}$ from $\{ t^i_+,
t^{i+1}_-\}$ and $\tilde{w}(q_i) =[a_i, b_i]^t$. By the assumption
that $x$ solves Thurston equation, $\prod_{i} a_i =\prod_{i} b_i$.
By the definition of $\phi_{\alpha_e}^{\delta_1}$, denoted by
$\psi$ for simplicity, we have $\psi(t_0^i)=\psi(t_0^{i+1}),
\psi(t_{\infty}^i)=\psi(t_{\infty}^{i+1})$,
$\psi(t^i_+)=\psi(t^i_-)$ and
\begin{equation}[\psi(t_0), \psi(t_{\infty}); \psi(t^i_+),
\psi(t^{i+1}_-)]=\tilde{w}(q_i). \end{equation}

%   \begin{figure}[ht!]
%\centering
%\includegraphics[scale=0.45]{fig2.eps}
%\caption{}
%\label{figure 2}
%\end{figure}

We claim that  $\psi(t_-^{m+1})$ defined by the identity (4) above is equal to
 $\psi(t_+^1)$. Indeed, by
corollary 3.3(5) and $\prod_{i} a_i =\prod_{i} b_i$ where all $a_i, b_i \in R^*$ and
the admissibility, we see that
 $[\psi(t_0), \psi(t_{\infty}); \psi(t^m_+), \psi(t^{m+1}_-)]
 =[  \psi(t_0), \psi(t_{\infty}); \psi(t^m_+), \psi(t^1_+)]$. By
 the uniqueness of the cross ratio, we conclude that
 $\psi(t_+^{m+1})=\psi(t_-^1)$.
\end{proof}

Now to define $\phi: \tilde{\Box} \to PR^1$, fix a tetrahedron $\sigma_0 \in \tilde{\mathcal T^*}$.
Let $t_1, t_2, t_3$ be three normal triangles in $\sigma_0$. Define $\phi(t_1)=[1,0]^t$,
$\phi(t_2)=[0,1]^t$ and $\phi(t_3)=[1,1]^t$ and use combinatorial continuation to
 define $\phi$  on $\tilde{\triangle}$.
 \end{proof}

 \subsection{The holonomy representation}

 Suppose $x$ is a solution to Thurston equation on $(M, \mathcal
 T)$ in a ring $R$ and $\phi:\tilde{\triangle} \to PR^1$ is an associated pseudo
 developing map. Then there exists a homomorphism $\rho:
 \pi_1(M^*) \to PSL(2,R)$    so that for all $\gamma \in
 \pi_1(M^*)$, considered as a deck transformation group for the
 universal cover $\pi_*: \tilde{M}^* \to M^*$,
 \begin{equation}
 \phi(\gamma) = \rho(\gamma) \phi.
 \end{equation}
 We call $\rho$ a \it holonomy representation \rm of $x$. It is
unique up to conjugation in $PSL(2,R)$. Here is the construction
of $\rho$. Fix an element $\gamma \in \pi_1(M^*)$. By the
construction, $\pi_1(M^*)$ acts on $\tilde{M}^*$, $\tilde{\mathcal
T}^*$, $\tilde{\triangle}$ and $\tilde{\Box}$ so that
$\pi_*(\gamma) =\pi_*$ for $\gamma \in \pi_1(M^*)$. This implies
$$[\phi(t_1), \phi(t_2);\phi(t_3), \phi(t_4)] =[\phi(\gamma t_1),
\phi(\gamma t_2); \phi(\gamma t_3), \phi(\gamma t_4)]$$ for all normal triangles $t_1, ..., t_4$ in
each tetrahedron $\sigma$ in $\tilde{T}^{(3)}$. By proposition 3.6, there exists an
element $\rho_{\sigma} (\gamma) \in PSL(2,R)$ so that
$$ \phi(\gamma t_i) =\rho_{\sigma}(\gamma) \phi(t_i)$$
where $t_i \subset \sigma$. We claim that $\rho_{\sigma}(\gamma) =
\rho_{\sigma'}(\gamma)$ for any two $\sigma, \sigma'
\in \tilde{\mathcal T^*}^{(3)}$. Indeed, since any two tetrahedra can be joint by an edge path in
the graph $G$, it suffices to show that $\rho_{\sigma}(\gamma) =
 \rho_{\sigma'}(\gamma)$ for two tetrahedra
sharing a codimension-1 face $\tau$. Let $t_1, t_2, t_3$ and $t_1', t_2', t_3'$ be the normal triangles
in $\sigma$ and $\sigma'$ respectively so that $t_i \sim t_i'$ and
 $t_1, t_2, t_3$ are adjacent to
$\tau$.   Now $\phi(t_i)=\phi(t_i')$ and $\gamma \phi(t_i) =\gamma
\phi(t_i')$, therefore, $\rho_{\sigma}(\gamma)
 \phi(t_i) =\rho_{\sigma'}(\gamma) \phi(t_i)$ for $i=1,2,3$. By the uniqueness part of
 proposition  3.6(2), it follows that $\rho_{\sigma}(\gamma)
 =\rho_{\sigma'}(\gamma)$. The common value is denoted by $\rho(\gamma)$.
 Given $\gamma_1, \gamma_2 \in
 \pi_1(M^*)$, by definition, $ \rho(\gamma_1 \gamma_2) \phi       =
  \phi(\gamma_1 \gamma_2)
 =\rho(\gamma_1) \phi(\gamma_2) =\rho(\gamma_1)\rho(\gamma_2)
 \phi$ and the uniqueness part of proposition 3.6, we see that
 $\rho(\gamma_1 \gamma_2) = \rho(\gamma_1) \rho(\gamma_2)$, i.e.,
 $\rho$ is a group homomorphism from $\pi_1(M^*)$ to $PSL(2,R)$.

 Note that the representation $\rho$ is trivial if and only if
 $\phi(\gamma t) = \phi(t)$ for all $t \in \tilde{\triangle}$ and $\gamma
 \in \pi_1(M^*)$. In this case, the pseudo developing map $\phi$ is defined on $\triangle
 \to PR^1$ so that $[\phi(t_1), \phi(t_2); \phi(t_3),
 \phi(t_4)]=[z(q), -z(q')]^t$. This was the construction
 in example 3.1.  In particular, the holomony representations
 associated to solutions in example 3.1 are trivial.

 \subsection{A proof of theorem 4.1}
 Suppose otherwise that there exists an edge $e \in \mathcal T^*$
 whose lift is an edge $e^*$ in $\tilde{\mathcal T}^*$ joining the
 same boundary component of $\tilde{M}^*$.   Take a tetrahedron
 $\sigma$ containing $e^*$ as an edge and let $t_1, t_2, t_3, t_4$
 be  all normal triangles in $\sigma$ so that $t_1, t_2$
 are adjacent to $e^*$. By definition, the pseudo developing map
 $\phi: \tilde{\triangle} \to PR^1$ satisfies the condition that  $\{ \phi(t_1),
 ..., \phi(t_4)\}$ is admissible. In particular, $\phi(t_1) \neq
 \phi(t_2)$. On the other hand, since $e^*$ ends at the same
 connected component of $\partial \tilde{M}^*$ which is a union of
 normal triangles related by sharing common normal arcs, there
 exists a sequence of normal triangles $s_1=t_1, s_2, ..., s_n
 =t_2$ in $\tilde{\triangle}$ so that $s_i \sim s_{i+1}$. In
 particular, $\phi(s_i) = \phi(s_{i+1})$. This implies that
 $\phi(t_1) = \phi(t_2)$ contradicting the assumption that
 $\phi(t_1) \neq \phi(t_2)$.

 To prove the second part of theorem 4.1 that $M$ is a closed
 connected 3-manifold, we first note that $\pi_1(M^*)$ is isomorphic
 to $\pi_1(M)$ under the homomorphism induced by inclusion. We will identify
 these two groups and identify
 $\tilde{M}^*$ as a $\pi_1(M)$ invariant subset of the universal cover $\tilde{M}$ of
 $M$. If $e$ is an edge in $\mathcal T$ ending at the same vertex
 $v$ in $\mathcal T$, let $\gamma \in \pi_1(M,v)$ be the deck
 transformation element corresponding to the loop $e$. We claim
 that $\rho(\gamma) \neq id$ in $PSL(2,R)$. Indeed, suppose
 $e^*$ is the lifting of $e$. Then by the statement just proved, $e^*$
 has two distinct vertices $u_1$ and $u_2$ in $\tilde{M}$ and $\phi(u_1) \neq
 \phi(u_2)$. By
 definition $\gamma(u_1) =u_2$. It follows that $\phi(u_2)
 =\phi(\gamma u_1) =\rho(\gamma) \phi(u_1)$. Since $\phi(u_1) \neq
 \phi(u_2)$, we obtain $\rho(\gamma) \neq id$.  This ends the
 proof.

\section{A universal construction}

Recall that $(M, \mathcal T)$ is a compact oriented triangulated
pseudo 3-manifold. The boundary $\partial M$ of $M$ is
triangulated by the subcomplex $\partial \mathcal T=\{ s \cap
\partial M | s \in \mathcal T\}$.  An edge in $\mathcal T$ is
called \it interior \rm if it is not in $\partial \mathcal T$. The
goal of this section is to introduce the \it Thurston ring \rm
$\mathcal R(\mathcal T)$ and its homogeneous version $\mathcal R_h
(\mathcal T)$. We will study the changes of $\mathcal R (\mathcal
T)$ when the triangulations are related by Pachner moves.

We will deal with quotients of the polynomial ring $\mathbf
Z[\Box]$ with $q \in \Box$ as variables. As a convention, we will
use $p \in \mathbf Z[\Box]$ to denote its image in the quotient
ring $\mathbf Z[\Box]/\mathcal I$.

\subsection{Thurston ring of a triangulation}
 The ``ground" ring in the
construction is the following. Let $\sigma$ be an oriented
tetrahedron  so that $q \to q' \to q''$ are the three quads in it.
Then the Thurston ring $\mathcal R (\sigma)$ is the quotient of
the polynomial ring $\mathbf Z[q,q',q'']$ modulo the ideal
generated by $q'(1-q)-1$, $q''(1-q')-1$, and $q(1-q'')-1$. Note
that this implies in $\mathcal R (\sigma)$, $q'=1/(1-q)$ and
$q''=(q-1)/q$ and furthermore, $\mathcal R(\sigma) \cong$ $\mathbf
Z[x, 1/x, 1/(1-x)]$ where $x$ is an independent variable.
Similarly, we defined $\mathcal R_h (\sigma)$ to be the quotient
ring $\mathbf Z[q,q',q'']/( q+q'+q'')$ where $(q+q'+q'')$ is the
ideal generated by $q+q'+q''$. Note that $\mathcal R_h (\sigma)
(\sigma) \cong \mathbf Z[x,y]$ the polynomial ring in two
independent variables.

 Recall that the tensor product
$R_1 \otimes R_2$ of two rings $R_1$ and $R_2$ is the tensor
product of  $R_1$ and $R_2$ considered as $\mathbf Z$ algebras.

\begin{defn} Suppose $(M, \mathcal T)$ is a compact oriented pseudo
3-manifold. The \it Thurston ring \rm $\mathcal R (\mathcal T)$ of
$\mathcal T$ is the quotient of the tensor product $\otimes_{
\sigma \in \mathcal T^{(3)}} \mathcal R(\sigma)$ modulo the ideal
generated by elements of the form $W_e-1$ where $W_e =\prod_{q
\sim e} q$ for all interior edges $e$. The \it homogeneous
Thurston ring \rm $\mathcal R_h (\mathcal T)$ is the quotient of
$\otimes_{\sigma \in \mathcal T^{(3)}} \mathcal R_h(\sigma)$
modulo the ideal generated by elements of the form $U_e = \prod_{
q \sim e} q -\prod_{ q \sim e} (-q')$, $q \to q'$, for all
interior edges $e$. The element $W_e=\prod_{q \sim e} q$  is
called the \it holonomy \rm of the edge $e$.
\end{defn}

By the construction, given a commutative ring $R$ with identity,
Thurston equation on $\mathcal T$ is solvable in $R$ if and only
if there exists a non-trivial ring homomorphism from $\mathcal
R(\mathcal T)$ to $R$. Therefore  theorem 1.1 can be stated as,
\begin{thm} Suppose $(M, \mathcal T)$ is a triangulated closed
connected 3-manifold so that one edge in $\mathcal T$ is a loop.
If $\mathcal R (\mathcal T) \neq \{0\}$, then $\pi_1(M) \neq
\{1\}$.
\end{thm}

 Note that $\mathcal R
(\mathcal T)$ is also the quotient $\mathbf Z[\Box]/\mathcal I$
where $\mathcal I$ is the ideal generated by $q'(1-q)-1$ for $q
\to q'$ and $W_e-1$ for interior edges $e$ and $\mathcal R_h
(\mathcal T)$ is the quotient of $\mathbf Z[\Box]/\mathcal I_h$
where $\mathcal I_h$ is the ideal generated by $q+q'+q''$ for $q
\to q' \to q''$ and $U_e$ for interior edges $e$.  We remark that
if there is $q \in \Box$ so that $q =0$ in $\mathcal R (\mathcal
T)$, then $\mathcal R(\mathcal T)=\{0\}$. Indeed, in this case the
identity element $1 = 1- q(1-q'')$ is in the ideal $\mathcal I$,
therefore, $\mathcal R(\mathcal T) =\{0\}$. In particular, if
$\mathcal T$ contains an interior edge of degree 1 (i.e., adjacent
to only one tetrahedron), then $\mathcal R(\mathcal T) =\{0\}$.

The relationship between $\mathcal R (\mathcal T)$ and $\mathcal
R_h (\mathcal T)$ is summarized in the following proposition. To
state it, recall that if $S$ is a multiplicatively closed subset
of a ring $R$, then the localization ring $R_S$ of $R$ at $S$ is
the quotient $R \times S/\sim$ where $(r_1, s_1) \sim (r_2, s_2)$
if there exists $s \in S$ so that $s(r_1s_2-r_2s_1)=0$. If $0 \in
S$, then $R_S=\{0\}$.

\begin{prop} Let $\mathcal S =\{ q_1....q_m |q_i \in \Box\}$ be
the multiplicatively closed subset of all monomials in $\Box$ in
$\mathcal R_h (\mathcal T)$. Then there exist a natural injective
ring homomorphism $F: \mathcal R (\mathcal T) \to \mathcal R_h
(\mathcal T)_S$ and a surjective ring homomorphism $G: \mathcal
R_h (\mathcal T) _S \to \mathcal R (\mathcal T)$ so that $GF=id$.
In particular, $\mathcal R (\mathcal T) =\{0\}$ if and only if
$\mathcal R_h (\mathcal T)_S=\{0\}$.
\end{prop}

\begin{proof} Define a ring homomorphism $F: \mathbf Z[\Box] \to
\mathcal R_h (\mathcal T)_S$ by $F(q) = -q/q'$ where $q \to q'$.
We claim that $F(\mathcal I) =\{0\}$ and thus $F$ induces a
homomorphism, still denoted by $F$, from $\mathcal R (\mathcal T)$
to $\mathcal R_h (\mathcal T)_S$. The generators of $\mathcal I$
are $q'(q-1)-1$ and $W_e-1$. If $q \to q' \to q''$, then
$F(q'(1-q)-1) = -\frac{q'}{q''}(1+\frac{q}{q'}) -1
=-\frac{q+q'+q''}{q''} =0$. For an interior edge $e$, $F(\prod_{q
\sim e} q -1 )= \frac{1}{\prod_{q \sim e} (-q')}$$(\prod_{q \sim
e} q - \prod_{q \sim e} (-q'))=0$. To construct the inverse of
$F$, we define $G: \mathbf Z[\Box] \to \mathcal R (\mathcal T)$ as
follows. For each tetrahedron $\sigma$ containing $q_1 \to q_2 \to
q_3$ where $q_1$ is specified, define $G(q_1) = q_1, G(q_2) = -1,
G(q_3) = 1-q_1$ in $\mathcal R (\mathcal T)$. By the construction,
we have $G(q) =-q G(q')$ for $q \to q'$ in $\Box$ and $\sum_{q
\subset \sigma} G(q)=0$ for each tetrahedron $\sigma$. We claim
that $G(\mathcal I_h) =\{0\}$, i.e., $G$ induces a ring
homomorphism, still denoted by $G: \mathcal R_h (\mathcal T) \to
\mathcal R (\mathcal T)$. Indeed, we have just verified the first
equation associated to each $\sigma \in \mathcal T$. For the
second type equation,  given any interior edge $e$, due to
$G(q)=-qG(q')$, $G(\prod_{ q \sim e} q -\prod_{q \sim e} (-q')) =
\prod_{q \sim e} G(q) -\prod_{q \sim e} (-G(q')) =[\prod_{q \sim
e} (-G(q'))] [\prod_{q \sim e} q -1] =0$. Note that by the
construction of $\mathcal R (\mathcal T)$, for $q \in \Box$, then
$q$ and $1-q$ are invertible in $\mathcal R (\mathcal T)$ with
inverses $1-q''$ and $q'$ where $q \to q' \to q''$. From the above
calculation, we see that $G$ induces a homomorphism from $\mathcal
R_h(\mathcal T)_S \to \mathcal R (\mathcal T)$. To check $GF=id$,
it suffices to see that $GF(q)=q$ for $q \in \Box$. Let $q \to q'
\to q''$ where $G(q')=-1$. Then $GF(q)=G(-q/q')=-q/(-1)=q$,
$GF(q')=G(-q'/q'')=1/(1-q)=q'$ and
$GF(q'')=G(-q''/q)=(q-1)/q=q''.$

To see the last statement, if $\mathcal R_h(\mathcal T)_S =\{0\}$,
then $\mathcal R(\mathcal T)=\{0\}$ since $F$ is injective. On the
other hand, if $\mathcal R(\mathcal T)=\{0\}$, then we claim that
$q=0$ in $\mathcal R_h(\mathcal T)$ for some $q \in \Box$. Indeed,
if not, then due to $F(q)=-q/q' \neq 0$, we see that $q \neq 0$ in
$\mathcal R(\mathcal T)$. This contradicts the assumption.
Therefore the multiplicatively closed set $S$ contains $0$. Hence
$\mathcal R_h(\mathcal T)_S =\{0\}$.

\end{proof}

%Note that if there is $q$ so that $q =0$ in  $\mathcal R (\mathcal
%T)$, then $\mathcal R(\mathcal T)=\{0\}$ since in this case the
%identity element $1 = 1- q(1-q'')$ is in the ideal $\mathcal I$.

%, \neq \{0\}$ if and only if $q=0$ for all
% $q \in \Box$. Therefore, if $\mathcal R (\mathcal T) = \{0\}$,
%then $\mathcal R_h (\mathcal T) =\mathbf Z$ by the proof above.

%Now if $\mathcal R_h (\mathcal T) \neq \{0\}$, there exists $q \in
%\Box$ so that $q \neq 0$ in $\mathcal R_h (\mathcal T)$. But
%$F(q)=-q/q'$ and $F$ is injective. Therefore $q \neq 0$ in
%$\mathcal R (\mathcal T)$.

\subsection{Pachner moves}
It is well known that any two triangulations of a closed pseudo
3-manifold are related by a sequence of Pachner moves \cite{Pa},
\cite{Pe}, \cite{Mat}. There are two types of Pachner moves: $1
\leftrightarrow 4$ move and $2 \leftrightarrow 3$ move. Two more
moves of types $0 \leftrightarrow 2_2$ and $0 \leftrightarrow 2_3$
are shown in figure 2. There moves create two new tetrahedra from
a triangle and a quadrilateral.  The $1 \leftrightarrow 4$ move is
a composition of a $0 \leftrightarrow 2_3$ move and a $2
\leftrightarrow 3$ move.

\begin{figure}[ht!]
\centering
\includegraphics[scale=0.6]{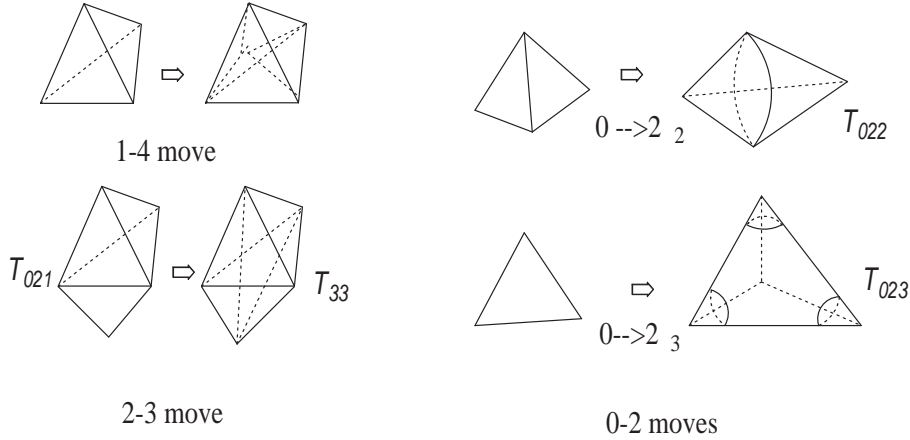}
\caption{Pachner moves} \label{figure 1}
\end{figure}

We will focus on the moves $0 \leftrightarrow 2_i$ $i=2,3$ and $2
\leftrightarrow 3$ in the rest of the paper and investigate their
effects on the Thurston ring $\mathcal R (\mathcal T)$. For this
purpose, we introduce the \it directed Pachner moves  \rm $i \to
j$ which means the Pachner move change a triangulation of fewer
tetrahedra to a triangulation of more tetrahedra, i.e., $0 \to
2_2$, $0 \to 2_3$ and $2 \to 3$. The following problem which
improves Pachner's theorem was investigated before. It dual
version for special spines was established by Makovetskii
\cite{Mak}. However, we are informed by S. Matveev \cite{Mat2}
that the following question is still open.

\medskip
\noindent {\bf Problem}. Suppose $\mathcal T_1$ and $\mathcal T_2$
are two triangulations of a closed pseudo 3-manifold $M$ so that
$\mathcal T_1^{(0)} =\mathcal T_2^{(0)}$. Then there exists a
third triangulation $\mathcal T$ of $M$ so that $\mathcal T$ is
obtained from both $\mathcal T_1$ and $\mathcal T_2$ by directed
Pachner moves $0 \to 2_2$ and $2 \to 3$.

\subsection{The directed moves $0 \to 2$ and $2 \to 3$}

There are four standard triangulations  of the 3-ball related to
the moves $0 \to 2_2$, $0\to 2_3$ and $2 \to 3$. Let $\mathcal
T_{021}, \mathcal T_{022}$ and $\mathcal T_{023}$ be the standard
triangulations of the 3-ball by two tetrahedra $\sigma^+$ and
$\sigma^-$. They are shown in figure 2 where the two tetrahedra in
$\mathcal T_{02i}$ share $i$ codimension-1 faces. The $2 \to 3$
move replaces $\mathcal T_{022}$ by $\mathcal T_{33}$. We will
calculate the ring $\mathcal R(\mathcal T_{02i})$ and $\mathcal
R(\mathcal T_{33})$ in this section. The results in subsection are
elementary and were known to experts in a less general setting.
See \cite{Se}, \cite{Ti} and others. We will emphasis the
naturality of the associated ring homomorphisms and holonomy
preserving properties.

For $\mathcal T_{02i}$, take a triangle in $\sigma^+ \cap
\sigma^-$ and let its edges be $e_1, e_2, e_3$ so that $e_3 \to
e_2 \to e_1$ in $\sigma^+$. Let $q_i^{\pm}$ be the quads in
$\sigma^{\pm}$ so that $q_i^{\pm} \sim e_i$.  Note that by the
construction  $q_1^+ \to q^+_2 \to q^+_3$ and $q_3^- \to q_2^- \to
q_1^-$.

\begin{lem} For $\mathcal T_{021}$, denote $q_1^+$ and $q^-_1$ by
$x,y$ in $\mathcal R(\mathcal T_{021})$ respectively. Then
$\mathcal R(\mathcal T_{021}) \cong \mathcal R(\sigma^+) \otimes
\mathcal R(\sigma^-)$ and the holonomies $W_{e_i}$ are: $W_{e_1} =
xy$, $W_{e_2} = \frac{y-1}{y-xy}$, $W_{e_3} =\frac{x-1}{x -xy}$.
The holonomies at all other edges are $x,y, 1/(1-x), 1/(1-y),
(x-1)/x, (y-1)/y$. In particular, if $W_{e_1}=1$, then $W_{e_i}
=1$ for $i=2,3$.
\end{lem}

\begin{proof} By definition, $q_2^+=1/(1-x), q_3^+=(x-1)/x$,
$q_2^-=(y-1)/y$ and $q^-_3 =1/(1-y)$. Since $W_{e_i} = q_i^+
q_i^-$, the result follows. \end{proof}

\begin{prop} (1) For the triangulation $\mathcal T_{023}$, the
inclusion homomorphism $\phi: \mathcal R(\sigma^{\pm}) \to
\mathcal R(\sigma^+) \otimes \mathcal R(\sigma^-)$ induces an
isomorphism $\Phi: \mathcal R(\sigma^{\pm}) \to \mathcal
R(\mathcal T_{023})$. Furthermore, the holonomy $W_e$ of each
boundary edge is $1$ in $\mathcal R(\mathcal T_{023})$.

(2) For the triangulation $\mathcal T_{022}$, let $e_0^{\pm} <
\sigma^{\pm}$ be the two boundary edges of degree 1 and assume
that $e_1$ is the interior edge. Then the inclusion homomorphism
$\phi: \mathcal R(\sigma^{\pm}) \to \mathcal R(\sigma^+) \otimes
\mathcal R(\sigma^-)$ induces an isomorphism $\Phi: \mathcal
R(\sigma^{\pm}) \to \mathcal R(\mathcal T_{022})$. Furthermore,
the holonomies $W_{e_0^{\pm}} = q_1^{\pm}$ and $W_e=1$ for all
other boundary edges $e$.

\end{prop}

\begin{figure}[ht!]
\centering
\includegraphics[scale=0.5]{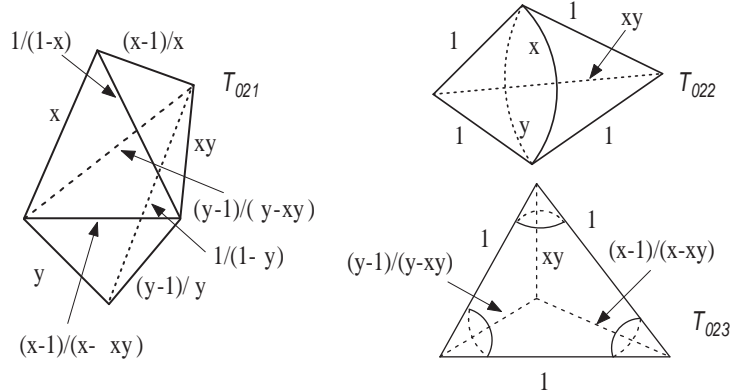}
\caption{Pachner moves and Thurston rings} \label{figure 1}
\end{figure}

\begin{proof}  To see part (1),  let $e_1$ and $e_2$ be the interior edges. Since
$W_{e_i} = q_i^+ q_i^-$, we have $q_i^-=1/q_i^+ \in \mathcal
R(\sigma^+)$ for $i=1,2$.  By lemma 5.4 and that $q_1^{\pm}
q_2^{\pm} q_2^{\pm}=-1$, we have $q_3^- =1/q_3^+$. It follows that
$\Phi: \mathcal R(\sigma^+) \to \mathcal R(\mathcal T_{023})$ is
onto. On the other hand, by exactly the same calculation as in
lemma 5.4, we see that $W_{e_2}=1$ and $W_{e_3}=1$ are consequence
of $W_{e_1}=1$, i.e., the ideal $\mathcal I$ in $\mathcal
R(\sigma^+) \otimes \mathcal R(\sigma^-)$ generated by $W_{e_1}-1$
contains $W_{e_2}-1$ and $W_{e_3}-1$. This shows that $\Phi$ is
injective. Therefore, $\Phi$ is a ring isomomorphism. Furthermore,
by definition, for each boundary edge $e^*_i$, the holonomy
$W_{e^*_i} = q_i^+ q_i^-=1$.

To see part (2), let $q^{\pm} = q_1^{\pm}$. By definition
$\mathcal R(\mathcal T_{022}) =\mathcal R(\sigma^+) \otimes
\mathcal R(\sigma^-)/(q^+ q^--1)$. In particular, $q^-=1/q^+$ and
that  $\Phi$ is an isomorphism.  The holonomies $W_e=1$ follow
from the definition and lemma 5.4.

\end{proof}

\begin{prop} The map $\phi: \Box(\mathcal T_{021}) \to \mathcal R(\mathcal T_{33})$ defined
by $\phi(\prod_{q \sim e} q) = \prod_{q \sim e} q$ for each
degree-1 edge $e$ induces a ring homomorphism $\Phi: \mathcal
R(\mathcal T_{021}) \to \mathcal R(\mathcal T_{33})$ so that for
all edges $e \in \mathcal T_{021}$, $\Phi(\prod_{q \sim e} q) =
\prod_{q \sim e} q$, i.e., $\Phi$ preserves holonomies.
Furthermore, let $S$ be the multiplicatively closed set consisting
of monomials in $x_iy_i-1$. Then $\Phi$ induces an isomorphism
from $\mathcal R(\mathcal T_{021})_S \to \mathcal R(\mathcal
T_{33})$.
\end{prop}
\begin{figure}[ht!]
\centering
\includegraphics[scale=0.45]{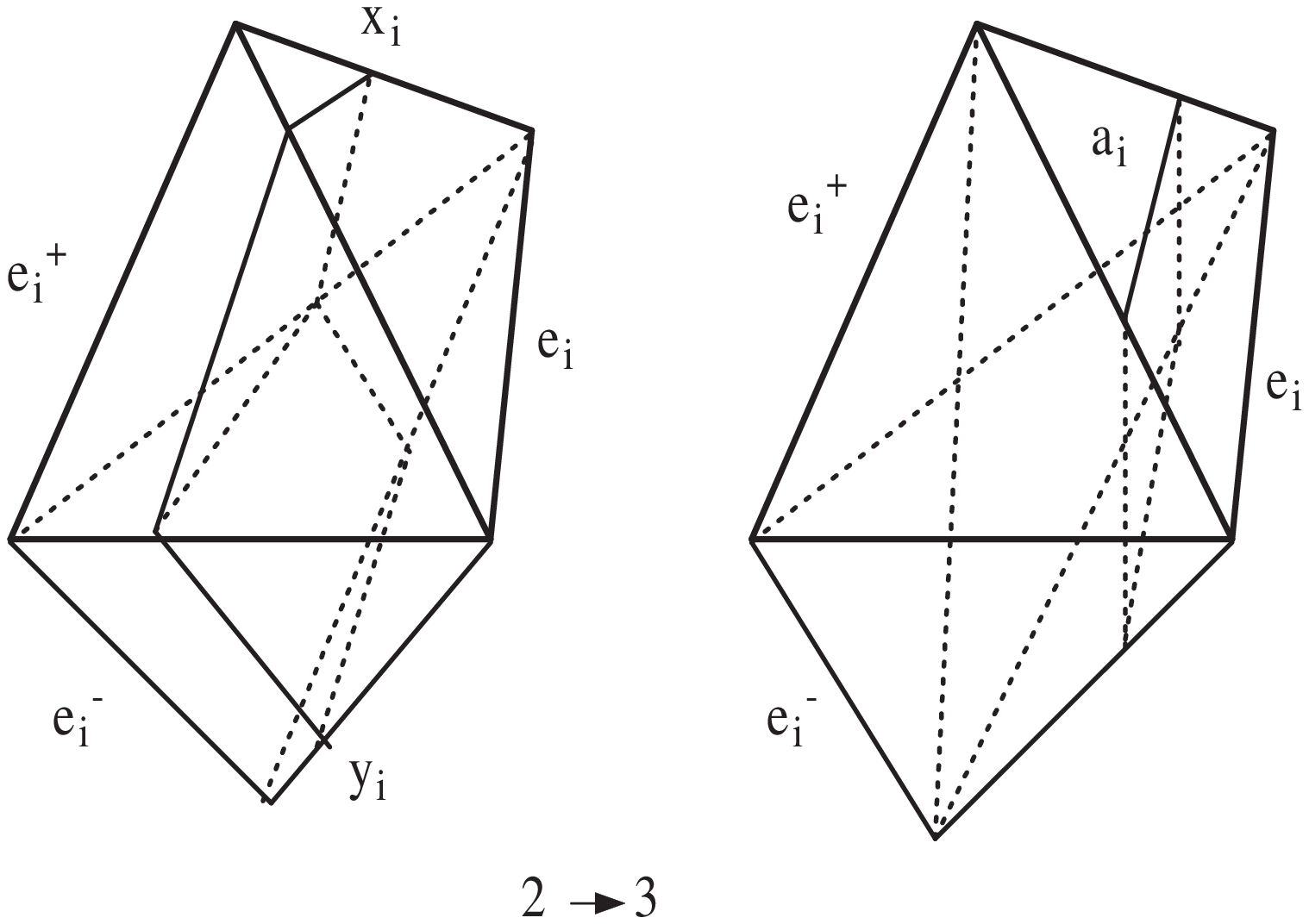}
\caption{$2 \to 3$ Pachner move} \label{figure 1}
\end{figure}

\begin{proof}  Let $e_i^{\pm}$ be the opposite edges of $e_i$ in
$\sigma^{\pm}$ and $\sigma_i$ be the tetrahedra in $\mathcal
T_{33}$ so that $e_i < \sigma_i$. The quads in $\sigma^{\pm}$ are
denoted by $x_i$ and $y_i$ so that $x_i \sim e_i$ and $y_i \sim
e_i$ respectively. The quads in $\sigma_i$ facing $e_i$ is denoted
by $a_i$. Let $b_i=a_i'$ and $c_i =b_i'$. Note that by the
construction $x_i'=x_{i+1}$ and $y'_{i+1}= y_i$. Furthermore, by
definition, $W_{e_i^+} (\mathcal T_{021}) = x_i$,
$W_{e_i^-}(\mathcal T_{021}) = y_i$, $W_{e_i}(\mathcal T_{021})=
x_iy_i$, $W_{e_i}(\mathcal T_{33}) = a_i$,  $W_{e_i^+}(\mathcal
T_{33}) = b_{i+2}c_{i+1}$, $W_{e_i^-}(\mathcal T_{33}) =
b_{i+1}c_{i+2}$, and $w_{e_0}(\mathcal T_{33})=a_1a_2a_3$ where
$e_0$ is the interior edge in $\mathcal T_{33}$. All indices are
calculated modulo 3.

By definition, the map $\phi: \Box(\mathcal T_{021}) \to \mathcal
T(\mathcal T_{33})$ is defined by $\phi(x_i) = b_{i+2}c_{i+1}$,
$\phi(y_i) =b_{i+1}c_{i+2}$. To show that $\phi$ induces a ring
homomorphism $\Phi: \mathcal R(\mathcal T_{021}) \to \mathcal
R(\mathcal T_{33})$, it suffices to show that
$\phi(x_{i+1})(1-\phi(x_i))=1$ and $\phi(y_i)(1-\phi(y_{i+1}) =1$.
Indeed, we have $b_i = 1/(1-a_i)$, $c_i=(a_i-1)/a_i$ and
$a_ia_{i+1}a_{i+2}=1$. Thus,
$$ \phi(x_{i+1})(1-\phi(x_i))= b_i c_{i+2}(1-b_{i+2}c_{i+1})$$
$$ =(\frac{1}{1-a_i} )(\frac{a_{i+2}-1}{a_{i+2}})
(1-(\frac{1}{1-a_{i+2}})( \frac{a_{i+1-1}}{a_{i+1}}))$$
$$=(\frac{a_{i+1} a_{i+2}}{a_{i+1} a_{i+2} -1})(
\frac{a_{i+2}-1}{a_{i+2}})( \frac{1
-a_{i+1}a_{i+2}}{(1-a_{i+2})a_{i+1}}) =1$$

The calculation for $\phi(y_i)(1-\phi(y_{i+1}) =1$ is very similar
and we omit the details.

To establish the identity $\phi(W_e) =W_e$ for the edges $e=e_i$,
i.e., we must verify that $\Phi(x_i y_i) = a_i$. Note that
$a_ib_ic_i=-1$ and $a_1a_2a_3=1$. Thus $\Phi(x_i y_i)
=b_{i+1}b_{i+2}c_{i+1}c_{i+2} =\frac{1}{a_{i+1} a_{i+2}} =a_i$.

To see that $\Phi$ induces an isomorphism from $\mathcal
R(\mathcal T_{021})_S$ to $\mathcal R(\mathcal T_{33})$, consider
the map $\psi: \Box(\mathcal T_{33}) \to \mathcal R(\mathcal
T_{021})$ so that $\psi(a_i) = x_iy_i$, $\psi(b_i)
=\frac{1}{1-x_iy_i}$ and $\psi(c_i)=\frac{x_iy_i-1}{x_iy_i}$. We
claim that $\psi$ induces a ring homomorphism $\Psi: \mathcal
R(\mathcal T_{021})_S \to \mathcal R(\mathcal T_{33})$ so that for
all edges $e$, $\Psi(\prod_{q \sim e} q) = \prod_{q \sim e} q$ and
$\Psi \Phi =id$, $\Phi \Psi=id$.

Indeed, to see that $\psi$ induces a ring homomorphism, we must
verify that $\Psi(a_1a_2a_3) =1$. This holds since
$\Psi(a_1a_2a_3) = x_1x_2x_3y_2y_2y_3 = (-1)(-1) =1$. To check
that $\Psi$ preserves the holonomies, it suffices to show $\Psi
(W_{e_i^{\pm}}) = W_{e_i^{\pm}}$. Since $x_{i+1} =
\frac{1}{1-x_i}$ and $y_i=\frac{1}{1-y_{i+1}}$, we have
$$\Psi(W_{e_i^+}) =\Psi(b_{i+2}c_{i+1})
=(\frac{1}{1-x_{i+2}y_{i+2}})(
\frac{x_{i+1}y_{i+1}-1}{x_{i+1}y_{i+1}})$$
$$=(\frac{1}{1-\frac{y_{i+1}-1}{(1-x_{i+1})y_{i+1}}})(
\frac{x_{i+1}y_{i+1}-1}{x_{i+1}y_{i+1}})$$
$$=\frac{(1-x_{i+1})y_{i+1}(x_{i+1}y_{i+1}-1)}{((y_{i+1}-x_{i+1}y_{i+1})
- y_{i+1}+1)x_{i+1}y_{i+1}}$$
$$=\frac{x_{i+1}-1}{x_{i+1}}
=x_i =W_{e_i^+}.$$

Essentially the same calculation shows $\Psi (W_{e_i^-}) =
W_{e_i^-}$.  Finally, due to holonomy preserving property of $\Phi
\Psi$ and $\Psi\Phi$, we have $\Psi \Phi =id$, $\Phi \Psi=id$.
\end{proof}

\subsection{Effects of Pachner moves}

Suppose $(M_i, \mathcal T_i)$ ($i=1,2$) are two compact
triangulated oriented pseudo 3-manifolds obtained as the quotients
$M_i = X_i/\sim_i$ of disjoint union $X_i$ of tetrahedra. Take $X
=X_1 \sqcup X_2$ and extend the identifications $\sim_i$ further
by identifying pairs of unidentified codimension-1 faces in $X$ by
orientation reversing affine homeomorphisms $\Phi$. The quotient
$X/\sim =M_1 \cup_{\Phi} M_2$ is called a gluing of $M_1$ and
$M_2$ along some subsurfaces of $\partial M_1$ and $\partial M_2$
by affine homeomorphism $\Phi$. The resulting triangulation will
be denoted by $\mathcal T_1 \cup_{\Phi} \mathcal T_2$. If $M_2
=\emptyset$, then $\Phi$ is a self-gluing of $M_1$. We denote the
result by $(M_1\cup_{\Phi}, \mathcal T_1\cup_{\Phi})$.  By
definition,
\begin{equation} \mathcal R(\mathcal T_1 \cup_{\Phi} \mathcal T_2)
=(\mathcal R(\mathcal T_1) \otimes \mathcal R(\mathcal
T_2))/\mathcal I \end{equation}
 where the ideal $\mathcal I$ is
generated by elements of the form $W_{e_1}(\mathcal T_1)
W_{e_2}(\mathcal T_2) -1$ with $e_1$ and $e_2$ being two boundary
edges which are identified to become an interior edge in $\mathcal
T_1 \cup_{\Phi} \mathcal T_2$. Note that there are natural ring
homomorphisms induced by the inclusion maps from $\Box(\mathcal
T_i)$ to $\Box(\mathcal T_1\cup_{\Phi} \mathcal T_2)$.

Using these notations, we can describe the effect of directed
Pachner moves $0 \to 2_3$, or $0 \to 2_2$ and $2 \to 3$ on
Thurston rings as follows. The moves $0 \to 2_3$ and $0 \to 2_2$
are of the form of replacing a self-glued $\mathcal T
\cup_{\Phi_1}$ by $\mathcal T \cup_{\Phi_2} \mathcal T_{02i}$ for
$i=2,3$. The move $2 \to 3$ replaces $\mathcal T \cup_{\Phi}
\mathcal T_{021}$ by $\mathcal T \cup_{\Phi} \mathcal T_{33}$.

Combining the definition (6) with the main results in \S5.3, we
have,

\begin{prop} Suppose $\mathcal T'$ is obtained from $\mathcal T$
by a directed Pachner move $0 \to 2_2$, $0 \to 2_3$ or $2 \to 3$.
Then there exists a holonomy preserving natural ring homomorphism
$\mathcal R(\mathcal T) \to \mathcal R(\mathcal T')$. \end{prop}

\section{Example of solving Thurston equation in finite rings}

Suppose $(M, \mathcal T)$ is a closed oriented pseudo 3-manifold
and $R$ is a commutative ring with identity and $x:\Box \to R$
solves Thurston equation. If $p$ is a prime number, let $F_{p^n}$
be the finite field of $p^n$ elements.

\begin{example} For $R =F_3$, then $x: \Box \to F_3-\{0,1\}$ is the
constant map $x(q)=2$. Thus, as mentioned in \S1.1, Thurston
equation is solvable if and only if each edge has even degree.
\end{example}

\begin{example} For $R=F_5=\{0,1,2,3,4\}$ and we are looking for
$x:\Box \to \{2,3,4\}$. Due to $1/(1-2)=4, 1/(1-4)=3, 4=2^2,
3=2^3$ so that $2^4=1$, we can write $x(q)= 2^{z(q)}$ where $z \in
\{1,2,3\}$. Thus Thurston equation is solvable if and only if for
$q \to q' \to q''$, $(z(q), z(q'), z(q'')) \in \{(1,2,3), (2,3,1),
(3,1,2) \}$ so that for each edge $e$, $\sum_{q \sim e} z(q) =0
\mod {4}$. \end{example}

\begin{example} For $R = F_{2^2} =\{0,1,a,b\}$ where $b=a+1=a^2$ and
$a^3=1$, we have $1/(1-a) =a$ and $1/(1-b)=b$. By writing solution
$x$ of Thurston equation as $x(q) = a^{z(q)}$ where $z(q) \in
\{1,2\}$, we see that $z(q) = z(q')$ if $q \to q'$. Therefore,
Thurston equation is solvable if and only if there is $z: \mathcal
T^{(3)} \to \{1,2\}$ so that for each edge $e$, $|\{ \sigma \in
\mathcal T^{(3)} | \sigma
>e, z(\sigma)=1\}| + 2|\{ \sigma \in \mathcal T^{(3)} | \sigma >e,
z(\sigma)=2\}| =0 \mod{3} $.
\end{example}

\begin{example}
For the field $F_7$, write $x(q)=3^{z(q)}$. Then Thurston equation
is solvable if and only if  $z: \Box \to \{1, 2,3,4,5\}$ satisfies
that $(z(q), z(q'), z(q'')) \in \{(1,1,1), (5,5,5), (2,3,4),
(3,4,2), (4,2,3)\}$ when $q \to q' \to q''$ and for each edge $e$,
$\sum_{q \sim e} z(q) =0 \mod{6}$.
\end{example}

\begin{example}
For the ring $\mathbf Z/9 \mathbf Z$ (not $F_{3^2}$), since a
solution $x(q)$ must satisfy $x(q)$ and $x(q)-1$ are invertible,
we conclude that $x(q) \in \{2,5,8\}=\{2, 2^5, 2^3\}$. Write
$x(q)=2^{z(q)}$. Therefore, Thurston equation is solvable if and
only if there is $z: \Box \to \{1,3,5\}$ so that  $(z(q), z(q'),
z(q'')) \in \{(1,3,5), (3,5,1), (5,1,3)\}$ if $q \to q' \to q''$
and for each edge $e$, $\sum_{q \sim e} z(q) =0 \mod{6}$. This
implies that the degree of each edge must be even.
 \end{example}

\begin{example}
For the ring $\mathbf Z/15 \mathbf Z$, the same argument as in
example 6.5 shows that Thurston equation is solvable if and only
if $x:\Box \to \{2,8, 14\}$ satisfies $(x(q), x(q'), x(q'')) \in
\{ (2, 14, 8), (14, 8, 2), (8, 2,14)\}$ if $q \to q' \to q''$ and
for each edge $e$, $\prod_{q \sim e} x(q) =0 \mod{15}$.
\end{example}

%\input sec8.tex

%%%%%%%%%%%%%%%%%%%%%%%%%%%%%%%%%%%%%%%%%%%%%%%%%%%%%%%%55

{}

\end{document}